\documentstyle[12pt,fleqn]{article}

\font\twlgot =eufm10 scaled \magstep1
\font\egtgot =eufm8
\font\sevgot =eufm7

\font\twlmsb =msbm10 scaled \magstep1
\font\egtmsb =msbm8
\font\sevmsb =msbm7

\newfam\gotfam
\def\pgot{\fam\gotfam\twlgot}
\textfont\gotfam\twlgot
\scriptfont\gotfam\egtgot
\scriptscriptfont\gotfam\sevgot
\def\got{\protect\pgot}

\newfam\msbfam
\textfont\msbfam\twlmsb
\scriptfont\msbfam\egtmsb
\scriptscriptfont\msbfam\sevmsb

\def\pBbb{\relax\ifmmode\expandafter\Bb\else\typeout{You cann't use
Bbb in text mode}\fi}
\def\Bb #1{{\fam\msbfam\relax#1}}

\newcommand{\gF}{{\got F}}

\def\thebibliography#1{\section*{References}\list
  {[\arabic{enumi}]}{\settowidth\labelwidth{#1}\leftmargin\labelwidth
    \advance\leftmargin\labelsep
    \usecounter{enumi}}
    \def\newblock{\hskip .11em plus .33em minus .07em}
    \sloppy\clubpenalty4000\widowpenalty4000
    \sfcode`\.=1000\relax}

\let\Large=\large

\def\op#1{\mathop{\fam0 #1}\limits}

\newcommand{\id}{{\rm Id\,}}

\newcommand{\di}{{\rm dim\,}}

\newcommand{\beq}{\begin{equation}}
\newcommand{\eeq}{\end{equation}}
\newcommand{\ben}{\begin{eqnarray}}
\newcommand{\een}{\end{eqnarray}}
\newcommand{\be}{\begin{eqnarray*}}
\newcommand{\ee}{\end{eqnarray*}}
\newcommand{\bea}{\begin{eqalph}}
\newcommand{\eea}{\end{eqalph}}

\newcommand{\cA}{{\cal A}}

\newcommand{\cT}{{\cal T}}

\newcommand{\cO}{{\cal O}}

\newcommand{\cF}{{\cal F}}

\newcommand{\cD}{{\cal D}}
\newcommand{\bL}{{\bf L}}
\newcommand{\bR}{{\bf R}}
\newcommand{\bC}{{\bf C}}

\newcommand{\bT}{{\bf T}}
\newcommand{\bZ}{{\bf Z}}

\newcommand{\fl}{\flat}
\newcommand{\sh}{\sharp}

\newcommand{\al}{\alpha}
\newcommand{\bt}{\beta}

\newcommand{\la}{\lambda}

\newcommand{\f}{\phi}

\newcommand{\Om}{\Omega}

\newcommand{\ve}{\varepsilon}
\newcommand{\g}{\gamma}
\newcommand{\G}{\Gamma}

\newcommand{\vt}{\vartheta}

\newcommand{\si}{\sigma}
\newcommand{\Si}{\Sigma}
\newcommand{\w}{\wedge}
\newcommand{\wt}{\widetilde}
\newcommand{\wh}{\widehat}
\newcommand{\ol}{\overline}
\newcommand{\dr}{\partial}
\newcommand{\ar}{\op\longrightarrow}
\newcommand{\ot}{\otimes}

\newcommand{\lng}{\langle}
\newcommand{\rng}{\rangle}

\newcounter{eqalph}
\newcounter{equationa}
\newcounter{theorem}
\newcounter{proposition}
\newcounter{lemma}
\newcounter{corollary}
\newcounter{definition}

\def\thedefinition{\arabic{definition}}

\newenvironment{proof}{\noindent {\it Proof.}}{\hfill $\Box$
\medskip }

\newenvironment{theo}{\refstepcounter{definition} \medskip\noindent
THEOREM \thedefinition.\it}{\medskip }
\newenvironment{prop}{\refstepcounter{definition} \medskip\noindent
PROPOSITION \thedefinition.\it}{\medskip }
\newenvironment{lem}{\refstepcounter{definition} \medskip\noindent  LEMMA
\thedefinition.\it }{\medskip }

\newenvironment{defi}{\refstepcounter{definition} \medskip\noindent 
DEFINITION \thedefinition.\it }{\medskip }

\newenvironment{eqalph}{\stepcounter{equation}
\setcounter{equationa}{\value{equation}}
\setcounter{equation}{0}

\begin{eqnarray}}{\end{eqnarray}\setcounter{equation}{\value{equationa}}}

\newcommand{\mar}[1]{}

\begin{document}
\hbox{}

{\parindent=0pt 

{ \Large \bf Geometric quantization of symplectic foliations}
\bigskip

{\sc G.SARDANASHVILY}

{ \small

{\it Department of Theoretical Physics,
Physics Faculty, Moscow State University, 117234 Moscow, Russia, e-mail:
sard@grav.phys.msu.su}
\bigskip

{\bf Abstract.}

Geometric quantization of a Poisson manifold need not imply
quantization of its symplectic leaves. We provide the leafwise
geometric quantization of a Poisson manifold, seen as a foliated one,
whose quantum algebra restricted to each leaf is quantized.

\medskip

\medskip


} }

\section{Introduction}

Though there is one-to-one correspondence between 
the (regular) Poisson structures on a smooth manifold and its
symplectic foliations, 
geometric quantization of a
Poisson manifold need not imply 
quantization of its symplectic leaves \cite{vais97}.

Firstly, contravariant connections fail to admit the pull-back operation.
Therefore, prequantization of a Poisson manifold does not determine
straightforwardly prequantization of its symplectic leaves.
Secondly, polarization of a Poisson manifold is
defined in terms of sheaves of functions, and it need not be associated 
to any distribution. As a consequence, its pull-back onto
a leaf is not polarization of a symplectic manifold in general.
Thirdly, a quantum algebra of a Poisson manifold contains 
the centre of a Poisson algebra. However, there are models 
where quantization of this centre has no physical meaning. For instance,
the centre of the Poisson algebra of a mechanical system with classical
parameters consists of functions of these parameters.

Geometric quantization of symplectic foliations disposes of 
these problems. The quantum algebra of a symplectic
foliation is also the quantum algebra of the associated Poisson
manifold such that its
restriction to each symplectic leaf is defined and quantized. 
Thus, geometric quantization of a
symplectic foliation provides the leafwise
quantization of a Poisson manifold. This is the case of systems whose
symplectic leaves are indexed by non-quantizable variables, e.g., systems
depending on classical parameters and constraint systems.

Geometric quantization of a symplectic foliation
is phrased in terms of the leafwise differential calculus
and leafwise connections on a foliated manifold. We show the  
following. Firstly, 
homomorphisms of the de Rham complex of a Poisson manifold both to its
Lichnerowicz--Poisson complex and the de Rham complex of its
symplectic leaf factorize through the leafwise de Rham complex, and
their cohomology groups do so. Secondly, any leafwise connection comes 
from a connection. Using these facts, we state the equivalence of
prequantization of a  
Poisson manifold to prequantization of its symplectic foliation, which also
yields prequantization of each symplectic leaf. 
On the contrary,
polarization of a symplectic foliation is associated to a particular
polarization of a Poisson manifold, and its restriction 
to any symplectic leaf is polarization of this leaf. 
Therefore, we define metaplectic correction of a symplectic foliation
so that its quantum algebra 
restricted to each
leaf is quantized. It is represented by Hermitian operators in the
pre-Hilbert space of leafwise half-forms, integrable over the leaves of this
foliation.

For example, the configuration space of 
a mechanical system with classical parameters is a 
fibre bundle $Q\to\Si$ over a
manifold of parameters 
$\Si$ \cite{book98,book00,sard00}. Its momentum phase space is the
vertical cotangent bundle $V^*Q$ of $Q\to\Si$ endowed with the
canonical Poisson structure, whose 
characteristic foliation coincides with the fibration $V^*Q\to\Si$.
One can think of its fibre $V^*_\si Q=T^*Q_\si$, $\si\in \Si$, as being the
momentum phase space of a mechanical system with fixed parameters. 
Of course, if a system is conservative, it can be quantized separately
at each point 
of a parameter space, but the leafwise quantization procedure need be called
into play if parameters are variable \cite{sard00}

\section{The leafwise differential calculus}

Manifolds throughout are assumed to be smooth, Hausdorff,
second-countable (i.e., paracompact), and connected.

A (regular) foliation $\cF$ of a manifold $Z$ consists
of (maximal) integral manifolds of an involutive distribution
$i_\cF:T\cF\to TZ$ on $Z$ \cite{rei}. A foliated manifold 
$(Z,\cF)$ admits an adapted coordinate atlas
\mar{spr850}\beq
\{(U_\xi;z^\la; z^i)\},\quad \la=1,\ldots,{\rm codim}\,\cF, \quad
i=1,\ldots,\di\cF,
\label{spr850}
\eeq
such that transition functions of coordinates $z^\la$ are independent of the 
remaining coordinates $z^i$ and,
for each leaf $F$ of a foliation $\cF$, the connected components
of $F\cap U_\xi$ are given by the equations
$z^\la=$const. These connected components and coordinates $(z^i)$ on
them make up a coordinate atlas of a leaf $F$.

The real Lie algebra $\cT_1(\cF)$ of 
global sections of the distribution $T\cF\to Z$ is
a $C^\infty(Z)$-submodule of the
derivation module of the $\bR$-ring
$C^\infty(Z)$ of smooth real functions on $Z$. Its kernel
$S_\cF(Z)\subset C^\infty(Z)$ consists of functions constant
on leaves of $\cF$. Therefore,
$\cT_1(\cF)$ is the Lie $S_\cF(Z)$-algebra of derivations of $C^\infty(Z)$,
regarded as a $S_\cF(Z)$-ring.
Then one can introduce the leafwise differential
calculus \cite{hect} as the 
Chevalley--Eilenberg differential calculus
over the $S_\cF(Z)$-ring $C^\infty(Z)$. It is defined as a subcomplex
\mar{spr892}\beq
0\to S_\cF(Z)\ar C^\infty(Z)\ar^{\wt d} \gF^1(Z) \cdots
\ar^{\wt d} \gF^{\di\cF}(Z) \to 0 \label{spr892}
\eeq
of the Chevalleqy--Eilenberg complex of
the Lie $S_\cF(Z)$-algebra $\cT_1(\cF)$
with coefficients in $C^\infty(Z)$ which consists of
$C^\infty(Z)$-multilinear skew-symmetric maps
$\op\times^r \cT_1(\cF) \to C^\infty(Z)$, $r=1,\ldots,\di\cF$
\cite{book00}. These maps are global sections of exterior products
$\op\w^r T\cF^*$ of the dual $T\cF^*\to Z$ of $T\cF\to Z$. 
They are called the leafwise forms
on a foliated manifold $(Z,\cF)$, and are given by the coordinate expression
\be
\f=\frac1{r!}\f_{i_1\ldots i_r}\wt dz^{i_1}\w\cdots\w \wt dz^{i_r},
\ee
where $\{\wt dz^i\}$ are the duals of the holonomic fibre
bases $\{\dr_i\}$ for $T\cF$. 
Then one can think of the 
Chevalley--Eilenberg coboundary operator 
\be
\wt d\f= \wt dz^k\w \dr_k\f=\frac{1}{r!}
\dr_k\f_{i_1\ldots i_r}\wt dz^k\w\wt dz^{i_1}\w\cdots\w\wt dz^{i_r}
\ee
as being the leafwise exterior differential.
Accordingly, (\ref{spr892}) is called the 
leafwise de Rham complex 
(or the tangential de Rham complex
in the terminology of \cite{hect}). This is the complex 
$(\cA^{0,*},d_f)$ in \cite{vais73}. Its cohomology $H^*_\cF(Z)$, called 
the leafwise de Rham cohomology, equals the cohomology 
$H^*(Z;S_\cF)$ of $Z$ with
coefficients in the sheaf $S_\cF$ of germs of elements of $S_\cF(Z)$
\cite{elk,most}. 

Let us consider the exact sequence
\mar{spr917}\beq
 0\to {\rm Ann}\,T\cF\ar T^*Z\ar^{i^*_\cF} T\cF^* \to 0. \label{spr917}
\eeq
of vector bundles over $Z$. Since it admits a splitting, 
the epimorphism $i^*_\cF$ yields an
epimorphism  of the graded algebra $\cO^*(Z)$ of exterior forms on $Z$ to
the algebra $\gF^*(Z)$ of leafwise forms. It obeys the condition 
$i^*_\cF\circ d=\wt d\circ i^*_\cF$ and, thereby, 
provides the cochain morphism
\mar{lmp04}\beq
 i^*_\cF: (\bR,\cO^*(Z),d)\to (S_\cF(Z),\cF^*(Z),\wt d), \qquad
 dz^\la\mapsto 0, \quad dz^i\mapsto\wt dz^i, \label{lmp04}
\eeq
of the de Rham complex of $Z$
to the leafwise de Rham complex (\ref{spr892}) and the 
corresponding homomorphism 
\mar{lmp00}\beq
[i^*_\cF]^*: H^*(Z)\to H^*_\cF(Z) \label{lmp00}
\eeq
of the de Rham cohomology of $Z$ to the leafwise one. Note that
$[i^*_\cF]^{r>0}$ need not be epimorphisms \cite{vais73}.

Given a leaf $i_F:F\to Z$ of a foliation $\cF$, we have the pull-back
homomorphism 
\mar{lmp30}\beq
(\bR,\cO^*(Z),d) \to (\bR,\cO^*(F),d) \label{lmp30}
\eeq
of the de Rham complex of $Z$ to that of $F$ and the corresponding
homomorphism of the de Rham cohomology groups
\mar{lmp31}\beq
H^*(Z) \to H^*(F). \label{lmp31}
\eeq

\begin{prop} \label{lmp32} \mar{lmp32} The homomorphisms (\ref{lmp30}) --
(\ref{lmp31}) factorize through the homomorphisms (\ref{lmp04}) --
(\ref{lmp00}).
\end{prop}

\begin{proof}
It is readily observed that
the pull-back bundles $i_F^* T\cF$ and $i_F^* T\cF^*$ over $F$ are
isomorphic to the tangent and the cotangent 
bundles of $F$, respectively. Moreover, a direct computation shows that
$i_F^*(\wt d\f)=d(i_F^*\f)$ for
any leafwise form $\f$. It follows that the cochain morphism (\ref{lmp30})
factorizes through the cochain morphism
(\ref{lmp04}) and the cochain morphism
\mar{lmp10}\beq
i^*_F: (S_\cF(Z),\cF^*(Z),\wt d)\to (\bR,\cO^*(F),d), \qquad
\wt dz^i\mapsto dz^i, \label{lmp10}
\eeq
of the leafwise de Rham complex of $(Z,\cF)$ to the de Rham complex
of $F$. Accordingly, the cohomology morphism (\ref{lmp31}) 
factorizes through the leafwise cohomology
\mar{lmp11}\beq
H^*(Z)\ar^{[i^*_\cF]} H^*_\cF(Z)\ar^{[i_F^*]} H^*(F). \label{lmp11}
\eeq
\end{proof}

Turn now to symplectic foliations.
Let $\cF$ be an even dimensional
foliation of a manifold $Z$.
A $\wt d$-closed non-degenerate leafwise
two-form $\Om_\cF$ on a foliated manifold $(Z,\cF)$ is called
symplectic. Its pull-back $i_F^*\Om_\cF$ onto
each leaf $F$ of $\cF$ is a 
symplectic form on $F$. 
If a leafwise symplectic form $\Om_\cF$ exists, it yields the
bundle isomorphism
\mar{spr896}\beq
\Om_\cF^\fl: T\cF\op\to_Z T\cF^*, \qquad
\Om_\cF^\fl:v\mapsto - v\rfloor\Om_\cF(z),
\qquad v\in T_z\cF.\label{spr896}
\eeq
The inverse isomorphism $\Om_\cF^\sh$ determines the
bivector field
\mar{spr904}\beq
w_\Om(\al,\bt)=\Om_\cF(\Om_\cF^\sh(i^*_\cF\al),\Om_\cF^\sh(i^*_\cF\bt)),
\qquad \forall \al,\bt\in T_z^*Z, \quad z\in Z,
\label{spr904}
\eeq
on $Z$ subordinate to $\op\w^2T\cF$. It is a
Poisson bivector field (see the relation 
(\ref{lmp12}) below).
The corresponding Poisson bracket reads
\mar{spr902}\beq
\{f,f'\}_\cF=\vt_f\rfloor \wt df', \qquad
\vt_f\rfloor\Om_\cF=-\wt df, \qquad \vt_f=\Om_\cF^\sh(\wt df),
\qquad f,f'\in C^\infty(Z).\label{spr902}
\eeq
Its kernel is $S_\cF(Z)$. 

Conversely, let $(Z,w)$ be a (regular) Poisson manifold and $\cF$ its
characteristic foliation. Since Ann$\,T\cF\subset T^*Z$ is
precisely the kernel of a Poisson bivector field $w$, 
the bundle homomorphism $w^\sh: T^*Z\op\to_Z TZ$
factorizes in a unique fashion 
\mar{lmp03}\beq
w^\sh:
T^*Z\ar_Z^{i^*_\cF} T\cF^*\ar_Z^{w^\sh_\cF}
T\cF\ar_Z^{i_\cF} TZ \label{lmp03}
\eeq
through the bundle isomorphism
\mar{lmp02}\beq
w_\cF^\sh: T\cF^*\op\to_Z T\cF,  \qquad
w^\sh_\cF:\al\mapsto -w(z)\lfloor \al, \qquad
\al\in T_z\cF^*. \label{lmp02}
\eeq
The inverse isomorphism $w_\cF^\fl$ yields the symplectic leafwise form 
\mar{spr903}\beq
\Om_\cF(v,v')=w(w_\cF^\fl(v),w_\cF^\fl(v')), \qquad \forall v,v'\in T_z\cF,
\qquad z\in Z. \label{spr903}
\eeq
The formulae (\ref{spr904}) and (\ref{spr903}) establish the 
above mentioned equivalence
between the Poisson structures on a manifold $Z$
and its symplectic foliations, though this equivalence need not be
preserved under morphisms. 

Let us consider the Lichnerowicz--Poisson
(henceforth LP) complex  
\mar{spr922}\beq
0\ar S_\cF(Z)\ar C^\infty(Z)\ar^{\wh w} \cT_1(Z)\ar^{\wh w} \cdots 
\ar^{\wh w}\cT_{\di\cF}(Z)
\to 0 \label{spr922}
\eeq
of multivector fields on a Poisson
manifold $(Z,w)$ with respect to the contravariant exterior differential
\be
\wh w: \cT_r(Z) \to \cT_{r+1}(Z), \qquad
\wh w(\vt)= -[w,\vt], \quad \vt\in \cT_*(Z), 
\ee
where $[.,.]$ denotes the Schouten--Nijenhuis bracket.
There are the cochain morphism
\mar{spr933}\ben
&& w^\sh: (\bR,\cO^*(Z),d) \to (S_\cF(Z),\cT_*(Z),-\wh w) \label{spr933}\\
&& w^\sh(\f)(\si_1,\ldots,\si_r)=
(-1)^r\f(w^\sh(\si_1),\ldots,w^\sh(\si_r)),
\quad \f\in\cO^r(Z),\quad \si_i\in \cO^1(Z), \nonumber\\
&&  \wh w\circ w^\sh= -w^\sh\circ d. \nonumber
\een
of the de Rham complex to the LP one and the corresponding 
homomorphism
\mar{gm111}\beq
[w^\sh]: H^*(Z)\to H^*_w(Z) \label{gm111}
\eeq
of the de Rham cohomology of $Z$ to the LP cohomology of the complex
(\ref{spr903}) \cite{vais}.

\begin{prop} \label{lmp15} \mar{lmp15}
The cochain morphism $w^\sh$ (\ref{spr933}) factorizes through the
leafwise complex (\ref{spr892}) and, accordingly, the cohomology
homomorphism $[w^\sh]$
(\ref{gm111}) does through the leafwise cohomology 
\mar{lmp60}\beq
H^*(Z)\ar^{[i^*_\cF]} H^*_\cF(Z)\ar H_w^*(Z). \label{lmp60}
\eeq
\end{prop}

\begin{proof}
Let $\cT_*(\cF)\subset \cT_*(Z)$ denote the exterior subalgebra of 
multivector fields
on $Z$ subordinate to $T\cF$, where $\cT_0(\cF)=C^\infty(Z)$. 
Clearly, $(S_\cF(Z),\cT_*(\cF),\wh w)$ is
a subcomplex of the LP complex (\ref{spr922}). 
Since 
\mar{lmp12}\beq
\wh w\circ \Om^\sh_\cF= -\Om^\sh_\cF\circ \wt d, \label{lmp12}
\eeq
the bundle isomorphism
$w^\sh_\cF=\Om^\sh_\cF$ 
(\ref{lmp02}) yields the cochain isomorphism
\be
\Om^\sh_\cF: (S_\cF(Z),\gF^*(Z),\wt d)\to (S_\cF(Z),\cT_*(\cF),-\wh w)
\ee
of the leafwise de Rham complex (\ref{spr892}) to the subcomplex
$(\cT_*(\cF),\wh w)$ of the LP complex (\ref{spr922}).
Then the composition
\mar{spr934}\beq
i_\cF\circ \Om^\sh_\cF: (S_\cF(Z),\gF^*(Z),\wt d)\to
(S_\cF(Z),\cT_*(Z),-\wh w) \label{spr934} 
\eeq
is a cochain monomorphism of the leafwise de Rham 
complex to the
LP one (\ref{spr922}). In view of the factorization (\ref{lmp03}), the
cochain morphism (\ref{spr933}) factorizes through  
the cochain morphisms (\ref{lmp04}) and (\ref{spr934}). Accordingly, the
cohomology homomorphism $[w^\sh]$ (\ref{gm111}) factorizes 
through the cohomology homomorphism $[i^*_\cF]$ (\ref{lmp00}) and the
cohomology homorphism
\mar{lmp05'}\beq
[i_\cF\circ \Om^\sh_\cF]: H^*_\cF(Z)\to H^*_w(Z). \label{lmp05'}
\eeq
\end{proof}

\section{Prequantization of a symplectic foliation}

Prequantization of a symplectic foliation $(\cF,\Om_\cF)$ of a manifold
$Z$ provides a representation 
\mar{lqm514}\beq
f\mapsto i\wh f, \qquad [\wh f,\wh f']=-i\wh{\{f,f'\}}_\cF, \label{lqm514}
\eeq
of the Poisson algebra $(C^\infty(Z),\{f,f'\}_\cF)$ by first order
differential operators on 
sections of a complex line 
bundle $\pi:C\to Z$. These operators
are given by the
Kostant--Souriau formula
\mar{lqq46}\beq
\wh f=-i\nabla_{\vt_f}^\cF +\ve f, \qquad \vt_f=\Om_\cF^\sh(\wt df),
\qquad \ve>0, \label{lqq46} 
\eeq
where $\nabla^\cF$ is a leafwise connection on $C\to Z$
such that its curvature form obeys the prequantization
condition
\mar{lmp61}\beq
\wt R=i\ve\Om_\cF. \label{lmp61}
\eeq
Using the fact that any leafwise connection comes from a connection
(see Theorem \ref{lmp42} below), we
will provide the cohomology analysis of this condition, and will show 
that prequantization of a symplectic foliation 
yields prequantization of its symplectic leaves. 

\begin{defi} \label{lmp100} \mar{lmp100}
In the framework of the leafwise differential calculus, 
a (linear) leafwise connection on the complex line bundle $C\to Z$ 
is defined as 
an algebraic connection $\nabla^\cF$ on the
$C^\infty(Z)$-module $C(Z)$ of global sections of this bundle, where
$C^\infty(Z)$ is regarded as a $S_\cF(Z)$-ring.
It associates to each element $\tau\in \cT_1(\cF)$ an $S_\cF(Z)$-linear
endomorphism $\nabla_\tau^\cF$ of $C(Z)$ which obeys the Leibniz rule
\mar{lmp55}\beq
\nabla_\tau^\cF(fs)=(\tau\rfloor\wt df)s +f\nabla_\tau^\cF(s), \qquad
f\in C^\infty(Z), \qquad s\in C(Z). \label{lmp55}
\eeq
\end{defi}

Recall that a linear connection on $C\to Z$
can equivalently be defined as an algebraic connection on the
module $C(Z)$ which assigns to each vector field $\tau\in
\cT_1(Z)$ on $Z$ an $\bR$-linear endomorphism of $C(Z)$
obeying the Leibniz rule (\ref{lmp55}). Restricted to $\cT_1(\cF)$, it
obviously yields a leafwise connection. In order to show that any
leafwise connection is of this form, we will
appeal to an alternative definition of a leafwise connection in terms
of leafwise forms.

The inverse
images $\pi^{-1}(F)$ of leaves $F$ of the foliation $\cF$ of $Z$
provide a (regular) foliation $C_\cF$ of the line bundle $C$. 
Given the (holomorphic)
tangent bundle $TC_\cF$ of this foliation, we have the
exact sequence of vector bundles
\mar{lmp18}\beq
0\to VC\ar_C TC_\cF\ar_C C\op\times_Z T\cF\to 0, \label{lmp18}
\eeq
where $VC$ is the (holomorphic) vertical tangent bundle of $C\to Z$.

\begin{defi} \label{lmp101} \mar{lmp101}
A (linear) leafwise 
connection on the complex line bundle $C\to Z$ is a splitting
of the exact sequence (\ref{lmp18}) which is linear over $C$.
\end{defi}

One can choose an adapted coordinate atlas $\{(U_\xi;z^\la; z^i)\}$
(\ref{spr850}) of a foliated manifold $(Z,\cF)$ such that $U_\xi$ are
trivialization domains of the complex line bundle $C\to Z$. Let 
$(z^\la; z^i;c)$, $c\in\bC$, 
be the
corresponding bundle coordinates on $C\to Z$. 
They are also adapted coordinates on the
foliated manifold $(C,C_\cF)$.  
With respect to these coordinates, a (linear) leafwise connection
is represented by a $TC_\cF$-valued leafwise
one-form
\mar{lmp21}\beq
A_\cF=\wt dz^i\ot(\dr_i +A_ic\dr_c), \label{lmp21}
\eeq
where $A_i$ are local complex functions on $C$. 

The exact sequence (\ref{lmp18}) is obviously a subsequence of the
exact sequence 
\be
0\to VC\ar_C TC\ar_C C\op\times_Z TZ\to 0,
\ee
where $TC$ is the holomorphic tangent bundle of $C$.
Consequently, any connection 
\mar{lmp103}\beq 
\G=dz^\la\ot(\dr_\la + \G_\la c\dr_c) + 
dz^i\ot(\dr_i +\G_ic\dr_c) \label{lmp103}
\eeq
on the complex line bundle $C\to Z$ yields a leafwise connection
\mar{lmp23}\beq
\G_\cF=\wt dz^i\ot(\dr_i +\G_ic\dr_c). \label{lmp23}
\eeq

\begin{theo} \label{lmp42} \mar{lmp42}
Any leafwise connection on the complex line bundle $C\to Z$ comes from 
a connection on it.
\end{theo}

\begin{proof}
Let $A_\cF$ (\ref{lmp21}) be a leafwise connection on 
$C\to Z$ and $\G_\cF$
(\ref{lmp23}) a leafwise connection which comes from some connection
$\G$ (\ref{lmp103}) on
$C\to Z$. Their affine difference over $C$ is a section
\be
Q=A_\cF-\G_\cF=\wt dz^i\ot(A_i -\G_i)c\dr_c
\ee
of the vector bundle $T\cF^*\op\ot_CVC\to C$.
Given some splitting 
\mar{lmp25}\beq
B: \wt dz^i \mapsto dz^i- B^i_\la dz^\la \label{lmp25}
\eeq
of the exact sequence (\ref{spr917}), the composition
\be
(B\ot \id_{VC})\circ Q=(dz^i- B^i_\la dz^\la)\ot(A_i -\G_i)c\dr_c: C\to
T^*Z\op\ot_C VC
\ee
is a soldering form on the complex line bundle $C\to Z$. Then 
\be
\G+(B\ot\id_{VC})\circ Q= 
dz^\la\ot(\dr_\la + [\G_\la -B^i_\la (A_i -\G_i)]c\dr_c) + 
dz^i\ot(\dr_i +A_ic\dr_c)
\ee
is a desired connection on $C\to Z$ which yields the leafwise 
connection $A_\cF$ (\ref{lmp21}).
\end{proof}

In particular, it follows that,
in view of the above
mentioned algebraic definition of a linear connection on a vector
bundle, Definition \ref{lmp100} and Definition \ref{lmp101} of
a leafwise connection are equivalent, namely, 
\be
\nabla^\cF s=\wt ds- A_i s\wt dz^i, \qquad s\in C(Z). 
\ee

The curvature of a leafwise connection $\nabla^\cF$ is defined as a
$C^\infty(Z)$-linear endomorphism
\mar{lmp07}\beq
\wt R(\tau,\tau')=\nabla_{[\tau,\tau']}^\cF- [\nabla_\tau^\cF,
\nabla_{\tau'}^\cF]=\tau^i
\tau'^j R_{ij}, \qquad
R_{ij}=\dr_i A_j-\dr_j A_i,
\label{lmp07}
\eeq
of $C(Z)$ for any vector fields $\tau,\tau'\in \cT_1(\cF)$.
It is represented by the
complex leafwise two-form
\mar{lmp08}\beq
\wt R=\frac12 R_{ij}\wt dz^i\w \wt dz^j. \label{lmp08}\\
\eeq
If a leafwise connection $\nabla^\cF$ comes from a connection $\nabla$,
its curvature leafwise form $\wt R$ (\ref{lmp08}) is the image 
$\wt R=i^*_\cF R$ 
of the curvature form $R$ of the connection $\nabla$ with respect to the
morphism $i^*_\cF$ (\ref{lmp04}).

Now let us turn to the prequantization condition (\ref{lmp61}). 

\begin{lem} \label{lmp63} \mar{lmp63}
Let us assume that there exists a leafwise connection
$\G_\cF$ on the complex line bundle $C\to Z$ which fulfils the
prequantization condition (\ref{lmp61}). Then,
for any Hermitian form $g$ on $C\to Z$, there exists
a leafwise connection $A_\cF^g$ on $C\to Z$ which: (i) satisfies the
condition
(\ref{lmp61}), (ii) preserves $g$, and (iii) comes from a
$U(1)$-principal connection on $C\to Z$.
\end{lem}

\begin{proof}
For any Hermitian form $g$ on $C\to Z$, there
exists an associated bundle atlas $\Psi^g=\{(z^\la;z^i,c)\}$ of $C$ with
$U(1)$-valued transition functions such that $g(c,c')=c\ol c'$.
Let the above mentioned leafwise connection $\G_\cF$ comes from 
a linear connection $\G$ (\ref{lmp103}) 
on $C\to Z$ written with respect to the atlas $\Psi^g$. 
The connection $\G$ is split into the sum $A^g + \g$ where 
\mar{lmp62}\beq
A^g=dz^\la\ot(\dr_\la + {\rm Im}(\G_\la)c\dr_c) + 
dz^i\ot(\dr_i +{\rm Im}(\G_i) c\dr_c) \label{lmp62}
\eeq
is a $U(1)$-principal connection, preserving the Hermitian form $g$. The
curvature forms $R$ of $\G$ and $R^g$ of $A^g$ 
obey the relation $R^g={\rm Im}(R)$. 
The connection $A^g$ (\ref{lmp62}) defines the leafwise
connection 
\mar{lmp73}\beq
A_\cF^g=i_\cF^*A= \wt dz^i\ot(\dr_i + iA^g_i c\dr_c), \qquad iA^g_i=
{\rm Im}(\G_i), \label{lmp73}
\eeq
preserving the Hermitian form $g$. Its curvature fulfils a desired
relation
\mar{lmp65}\beq
\wt R^g=i_\cF^*R^g={\rm Im}(i_\cF^*R)= i\ve\Om_\cF. \label{lmp65}
\eeq
\end{proof}

Since $A^g$ (\ref{lmp62}) is a $U(1)$-principal connection, its
curvature form $R^g$ is related to the first Chern form 
of integer de Rham cohomology class
by the formula
$c_1=i(2\pi)^{-1}R^g$.
If the prequantization
condition (\ref{lmp61}) holds,
the relation (\ref{lmp65}) shows that the leafwise
cohomology class of the leafwise 
form $(2\pi)^{-1}\ve\Om_\cF$ is the image of an integer de Rham
cohomology class with respect to the cohomology morphism $[i^*_\cF]$
(\ref{lmp00}). Conversely, if a leafwise
symplectic form $\Om_\cF$ on a foliated manifold $(Z,\cF)$ is of this
type, there exists a complex 
line bundle $C\to Z$ and a $U(1)$-principal connection $A$ on $C\to Z$
such that the leafwise connection $i^*_\cF A$ fulfils the relation
(\ref{lmp61}). Thus, we have stated the following.

\begin{prop} \label{lmp66} \mar{lmp66}
A symplectic foliation $(\cF,\Om_\cF)$ of a manifold $Z$ admits
prequantization (\ref{lqq46}) iff  
the leafwise cohomology class of $(2\pi)^{-1}\ve\Om_\cF$ is the image
of an integer de Rham 
cohomology class of $Z$.
\end{prop}

In particular, let $(Z,w)$ be a Poisson manifold and $(\cF,\Om_\cF)$
its characteristic symplectic foliation. As is well-known, a Poisson manifold
admits prequantization iff the LP cohomology class of the
bivector field $(2\pi)^{-1}\ve w$, $\ve>0$, is the image of an
integer de Rham cohomology class 
with respect to the cohomology morphism $[w^\sh]$ (\ref{gm111})
\cite{vais91,vais}. By virtue of Proposition \ref{lmp15}, this morphism
factorizes through the cohomology morphism $[i^*_\cF]$
(\ref{lmp00}). Therefore, in accordance with Proposition \ref{lmp66},
prequantization of a Poisson manifold takes place iff prequantization
of its symplectic foliation does well, and both these prequantizations 
utilize the same
prequantization bundle $C\to Z$. Herewith, each leafwise connection
$\nabla^\cF$ obeying the prequantization condition (\ref{lmp61})
yields the admissible contravariant connection
$\nabla_\f^w=\nabla_{w^\sh(\f)}^\cF$, $\f\in \cO^1(Z)$, 
on $C\to Z$ whose curvature bivector equals $i\ve w$. Clearly,
$\nabla^\cF$ and $\nabla^w$ lead to the same
prequantization formula (\ref{lqq46}).

Let $F$ be a leaf of a symplectic foliation
$(\cF,\Om_\cF)$ provided with the symplectic form $\Om_F=i^*_F\Om_\cF$.
In accordance with Proposition \ref{lmp32} and
the commutative diagram of cohomology groups
\be
\begin{array}{ccc}
 H^*(Z;\bZ) &\ar & H^*(Z)\\
 \put(0,10){\vector(0,-1){20}} & & \put(0,10){\vector(0,-1){20}}\\
H^*(F;\bZ) &\ar & H^*(F)
\end{array},
\ee
the symplectic form $(2\pi)^{-1}\ve \Om_F$ belongs to an integer
de Rham cohomology 
class if a leafwise symplectic form $\Om_\cF$ fulfils the condition of
Proposition \ref{lmp66}. This states the following.

\begin{prop} \label{lmp70} \mar{lmp70}
If a symplectic foliation 
admits prequantization, each its symplectic leaf
does well.
\end{prop}

The corresponding prequantization bundle for $F$ is the pull-back
complex line bundle $i^*_FC$, coordinated by $(z^i,c)$. Furthermore,
let $A_\cF^g$ (\ref{lmp73}) be 
a leafwise connection on the prequantization bundle $C\to Z$
which obeys Lemma \ref{lmp63}, i.e., comes from a
$U(1)$-principal connection $A^g$ on $C\to Z$. Then the pull-back 
\mar{lmp130}\beq
A_F=i^*_FA^g=dz^i\ot(\dr_i +ii^*_F(A^g_i)c\dr_c) \label{lmp130}
\eeq
of the connection $A^g$ onto 
$i^*_FC\to F$ satisfies the 
prequantization condition
\be
R_F=i^*_FR=i\ve \Om_F,
\ee
and preserves the pull-back Hermitian form $i^*_Fg$ on $i^*_\cF C\to
F$. 

\section{Polarization of a symplectic foliation}

Let us define polarization of a symplectic foliation $(\cF,\Om_\cF)$ of
a manifold $Z$ as a maximal (regular) involutive distribution
$\bT\subset T\cF$ on $Z$ such that
\mar{lmp71}\beq
\Om_\cF(u,v)=0, \qquad \forall u,v\in\bT_z, \qquad z\in Z. \label{lmp71}
\eeq
Given  the Lie algebra $\bT(Z)$ of $\bT$-subordinate vector fields on
$Z$, let $\cA_\cF\subset C^\infty(Z)$ be the complexified 
subalgebra of functions $f$
whose leafwise Hamiltonian vector fields $\vt_f$ (\ref{spr902}) fulfil
the condition $[\vt_f,\bT(Z)]\subset \bT(Z)$.
It is called the quantum algebra of a symplectic
foliation $(\cF,\Om_\cF)$ with respect to the polarization $\bT$. 
This algebra obviously contains the centre $S_\cF(Z)$
of the Poisson algebra $(C^\infty(Z),\{,\}_\cF)$, and is a Lie
$S_\cF(Z)$-algebra. 

\begin{prop} \label{lmp72} \mar{lmp72}
Every polarization $\bT$ of a symplectic foliation $(\cF,\Om_\cF)$
yields polarization of the associated Poisson manifold $(Z,w_\Om)$.
\end{prop}

\begin{proof} 
Let us consider the presheaf of local smooth functions $f$ 
on $Z$ whose leafwise Hamiltonian 
vector fields $\vt_f$ (\ref{spr902}) are subordinate to $\bT$. 
The sheaf $\Phi$ of germs of these functions is 
polarization of the Poisson manifold $(Z,w_\Om)$. 
Equivalently, $\Phi$ is the sheaf of germs of functions on $Z$ whose
leafwise differentials are subordinate to the codistribution 
$\Om_\cF^\fl\bT$. 
\end{proof}

Note that the polarization $\Phi$ 
need not be maximal, unless $\bT$ is of maximal dimension $\di\cF/2$.
It belongs to the following particular type of polarizations of a
Poisson manifold. Since the cochain morphism $i^*_\cF$ (\ref{lmp04}) is
an epimorphism, the leafwise differential calculus $\gF^*$ is
universal, i.e., the leafwise differentials $\wt df$ of functions $f\in
C^\infty(Z)$ on
$Z$ make up a basis for the $C^\infty(Z)$-module $\gF^1(Z)$. Let 
$\Phi(Z)$ denote the structure $\bR$-module of global sections of the
sheaf $\Phi$. Then the leafwise differentials of elements of $\Phi(Z)$
make up a basis for the $C^\infty(Z)$-module of global sections of the
codistribution $\Om_\cF^\fl\bT$. Equivalently, the leafwise Hamiltonian
vector fields of elements of $\Phi(Z)$ constitute a basis for the 
$C^\infty(Z)$-module $\bT(Z)$. Then one can easily show that
polarization $\bT$ of a symplectic foliation $(\cF,\Om_\cF)$ and the
corresponding polarization $\Phi$ of the Poisson manifold $(Z,w_\Om)$
in Proposition \ref{lmp72} define the same quantum algebra $\cA_\cF$.

Let $(F,\Om_F)$ be a symplectic leaf of a symplectic foliation
$(\cF,\Om_\cF)$. Given a polarization $\bT\to Z$ of $(\cF,\Om_\cF)$, its
restriction $\bT_F=i^*_F\bT\subset i^*_FT\cF=TF$ to $F$ is an
involutive distribution on $F$. It obeys the condition
\be
i^*_F\Om_\cF(u,v)=0, \qquad \forall u,v\in\bT_{Fz}, \qquad z\in F,
\ee
i.e., is polarization of the symplectic manifold $(F,\Om_F)$.
Thus, we have stated the following.

\begin{prop} \label{lmp75} \mar{lmp75}
Polarization of a symplectic foliation defines polarization of each
symplectic leaf. 
\end{prop}

Clearly, the quantum algebra $\cA_F$ of a symplectic leaf $F$ with respect to
the polarization $\bT_F$ contains all elements $i^*_Ff$ of the quantum algebra
$\cA_\cF$ restricted to $F$.

\section{Quantization of a symplectic foliation}

Since $\cA_\cF$ is the quantum algebra both of a symplectic foliation
$(\cF,\Om_\cF)$ and the associated Poisson manifold $(Z,w_\Om)$, 
one let us start from the
standard metaplectic correction technique \cite{eche98,wood}. 

Assuming that $Z$ is oriented and that 
$H^2(Z;\bZ_2)=0$, let us consider the metalinear complex
line bundle 
$\cD\to Z$ characterized by an atlas
$\Psi_Z=\{(U_\xi;z^\la;z^i;c)\}$ with 
the transition functions $c'=Sc$ such that $S\ol 
S$ is the inverse Jacobian of coordinate transition functions on $Z$.
Global sections  
of this bundle are half-forms on $Z$. The metalinear bundle
$\cD$ belongs to the category of natural bundles, and the Lie derivative
\mar{lmp78}\beq
\bL_\tau=\tau^\la \dr_\la + \tau^i\dr_i+\frac12(\dr_\la\tau^\la+
\dr_i\tau^i) \label{lmp78}
\eeq
of its sections along any vector field $\tau$ on $Z$ is defined.
The quantization bundle is the tensor product $Y=C\ot\cD$. 
The space $Y_K(Z)$ of
its sections of  
compact support is provided
with the non-degenerate Hermitian form 
\mar{lmp77}\beq
\lng\rho|\rho'\rng=\left(\frac1{2\pi}\right)^{\di Z/2}\op\int_Z
\rho\rho', \qquad \rho,\rho'\in Y_K(Z),\label{lmp77}
\eeq
written with respect to the atlases $\Psi^g$ of $C$ and $\Psi_Z$ of
$\cD$. Given the leafwise connection $A^g_\cF$ (\ref{lmp73})
and the Lie derivative
$\bL$ (\ref{lmp78}), one can assign the first order differential operator  
\mar{lmp80}\beq
\wh f=-i[(\nabla_{\vt_f}^\cF +i\ve f)\ot\id +\id\ot\bL_{\vt_f}]=
-i[\nabla_{\vt_f}^\cF +i\ve f+\frac12\dr_i\vt_f^i], \qquad f\in\cA_\cF,
\label{lmp80} 
\eeq
on $Y_K(Z)$ to
each element of the quantum algebra $\cA_\cF$. These
operators obey the Dirac condition (\ref{lqm514}), and 
provide a representation of the quantum algebra $\cA_\cF$ by
(unbounded) Hermitian operators in the pre-Hilbert space $Y_K(Z)$.
Finally, this representation is restricted to the subspace
$E$ of sections $\rho\in Y_K(Z)$ which obey the condition  
\be
(\nabla_\vt^\cF\ot\id +\id\ot\bL_\vt)\rho=(\nabla_\vt^\cF 
+\frac12\dr_i\vt^i)\rho=0
\ee
for all $\bT$-subordinate leafwise Hamiltonian
vector fields $\vt$.

However, it may happen that 
the above quantization has no physical sense because 
the Hermitian form
(\ref{lmp77}) on the carrier space $E$ and, consequently, the mean
values of operators (\ref{lmp80}) are defined by 
integration over the whole manifold $Z$. For instance, it implies
integration over classical parameters. Therefore, we
suggest a different scheme of quantization of symplectic foliations. 

Let us consider the exterior bundle $\op\w^mT\cF^*$, $m=\di\cF$.
Its structure group $GL(m,\bR)$
is reducible to the group $GL^+(m,\bR)$ since a symplectic
foliation is oriented. One can regard this fibre bundle as being
associated to a $GL(m,\bC)$-principal bundle $P\to Z$. 
As earlier, let us assume that $H^2(Z;\bZ_2)=0$. Then the principal
bundle $P$ admits a two-fold covering principal bundle with the 
structure metalinear group $ML(m,\bC)$ \cite{eche98}. 
As a consequence, there
exists a complex line bundle $\cD_\cF\to Z$ 
characterized by an atlas $\Psi_\cF=\{(U_\xi;z^\la;z^i;c)\}$ with
the transition functions $c'=S_\cF c$ such that 
\be
S_\cF\ol S_\cF=\det\left(\frac{\dr z^i}{\dr z'^j}\right).
\ee
One can think of its sections as being leafwise half-forms on $Z$.
The metalinear bundle $\cD_\cF\to Z$ admits the canonical lift 
of any $\bT$-subordinate vector field $\tau$ on $Z$. 
The corresponding
Lie derivative of its sections reads
\mar{lmp82}\beq
\bL_\tau^\cF=\tau^i\dr_i+\frac12\dr_i\tau^i. \label{lmp82}
\eeq

We define the quantization bundle as the tensor product
$Y_\cF=C\ot\cD_\cF$.
Given a leafwise connection $A^g_\cF$ (\ref{lmp73})
and the Lie derivative $\bL^\cF$ (\ref{lmp82}), let us 
associate  the first order differential operator 
\mar{lmp84}\beq
\wh f=-i[(\nabla_{\vt_f}^\cF +i\ve f)\ot\id +\id\ot\bL_{\vt_f}^\cF]=
-i[\nabla_{\vt_f}^\cF +i\ve f+\frac12\dr_i\vt_f^i], \qquad f\in\cA_\cF,
\label{lmp84}
\eeq
on sections $\rho_\cF$ of $Y_\cF$ to
each element of the quantum algebra $\cA_\cF$. A direct
computation with respect to the local Darboux coordinates on $Z$ proves
the following. 

\begin{lem} \label{lmp120} \mar{lmp120}
The operators (\ref{lmp84}) obey the Dirac condition (\ref{lqm514}).
\end{lem}

\begin{lem} \label{lmp121} \mar{lmp121}
If a section $\rho_\cF$ fulfils the condition
\mar{lmp86}\beq
(\nabla_\vt^\cF\ot\id +\id\ot\bL_\vt^\cF)\rho_\cF=(\nabla_\vt^\cF 
+\frac12\dr_i\vt^i)\rho_\cF=0  \label{lmp86}
\eeq
for all $\bT$-subordinate leafwise Hamiltonian vector field $\vt$, 
then $\wh f\rho_\cF$ for any $f\in\cA_\cF$ possesses the same property.
\end{lem}

Let us restrict the representation of the 
quantum algebra $\cA_\cF$ by the operators (\ref{lmp84}) to the subspace 
$E_\cF\in Y_\cF(Z)$ of
sections $\rho_\cF$ which obey the condition (\ref{lmp86}) and whose
restriction to 
any leaf of $\cF$ is of compact support. The last condition is
motivated by the following.

Since $i^*_FT\cF^*=T^*F$, the pull-back $i^*_F\cD_\cF$
of $\cD_\cF$ onto a leaf $F$ is a metalinear
bundle of half-forms on $F$. By virtue of
Proposition \ref{lmp70} and Proposition \ref{lmp75}, the pull-back
$i^*_FY_\cF$ of the quantization bundle $Y_\cF\to Z$ onto
$F$ is a quantization bundle for the symplectic manifold
$(F,i^*_F\Om_\cF)$. Given the pull-back connection $A_F$ (\ref{lmp130})
and the polarization $\bT_F=i^*_F\bT$, this symplectic manifold is subject to
the standard geometric quantization by the
first order differential operators 
\mar{lmp92}\beq
\wh f=-i(i_F^*\nabla_{\vt_f}^\cF +i\ve f +\frac12\dr_i\vt_f^i), \qquad
f\in \cA_F,  \label{lmp92}
\eeq
on sections $\rho_F$ of $i^*_FY_\cF\to F$ of compact support 
which obey the condition 
\mar{lmp133}\beq
(i_F^*\nabla_\vt^\cF +\frac12\dr_i\vt^i)\rho_F=0 \label{lmp133}
\eeq
for all $\bT_F$-subordinate Hamiltonian vector fields $\vt$ on $F$.
These sections constitute a pre-Hilbert space $E_F$ with respect to 
the Hermitian form
\be
\lng\rho_F|\rho'_F\rng=\left(\frac1{2\pi}\right)^{m/2}\op\int_F
\rho_F \rho'_F.
\ee
The key point is the following.

\begin{prop} \label{lmp132} \mar{lmp132}
We have $i^*_FE_\cF\subset E_F$, and the relation 
\beq
i^*_F(\wh f\rho_\cF)=\wh{(i^*_Ff)}(i^*_F\rho_\cF) \label{lmp93}
\eeq
holds for all elements $f\in\cA_\cF$ and $\rho_\cF\in E_\cF$.
\end{prop}

\begin{proof} One can use the fact that the expressions (\ref{lmp92})
and (\ref{lmp133}) have the 
same coordinate form as the expressions (\ref{lmp84}) and
(\ref{lmp86}) where $z^\la=$const.
\end{proof}

The relation (\ref{lmp93}) enables one to think of the operators $\wh f$
(\ref{lmp84}) in $E_\cF$ as being the leafwise quantization of the
$S_\cF(Z)$-algebra $\cA_\cF$ in the pre-Hilbert $S_\cF(Z)$-module. 

For example, the instantwise
quantization of time-dependent mechanics is of this type 
\cite{gi01,sni}.

\end{document}